\documentclass[11pt]{article}
\usepackage[usenames]{color}
\usepackage{amssymb,psfig,epsfig,latexsym,graphicx}
\usepackage{amsmath}

\usepackage[english]{babel}
\usepackage{setspace}
\usepackage{dsfont}

\definecolor{webgreen}{rgb}{0,.5,0}
\definecolor{webbrown}{rgb}{.6,0,0}

\newtheorem{theorem}{Theorem}

\newtheorem{prop}[theorem]{Proposition}

\newcommand{\eqn}[1]{(\ref{#1})}
\newcommand{\bsq}{{\vrule height .9ex width .8ex depth -.1ex }}

\newcommand{\La}{\Lambda}

\newcommand{\eps}{{\epsilon}}

\newcommand{\si}{\sigma}
\newcommand{\eeq}{\end{equation}}

\newcommand{\beql}[1]{\begin{equation}\label{#1}}
\DeclareMathOperator{\tr}{tr}
\newcommand{\RR}{\mathbb R}
\newcommand{\QQ}{\mathbb Q}
\newcommand{\ZZ}{\mathbb Z}

\makeatletter
\def\@sect#1#2#3#4#5#6[#7]#8{\ifnum #2>\c@secnumdepth
 \def\@svsec{}\else
 \refstepcounter{#1}\edef\@svsec{\csname the#1\endcsname.\hskip .75em }\fi
 \@tempskipa #5\relax
 \ifdim \@tempskipa>\z@
 \begingroup #6\relax
 \@hangfrom{\hskip #3\relax\@svsec}{\interlinepenalty \@M #8\par}%
 \endgroup
 \csname #1mark\endcsname{#7}\addcontentsline
 {toc}{#1}{\ifnum #2>\c@secnumdepth \else
 \protect\numberline{\csname the#1\endcsname}\fi
 #7}\else
 \def\@svsechd{#6\hskip #3\@svsec #8\csname #1mark\endcsname
 {#7}\addcontentsline
 {toc}{#1}{\ifnum #2>\c@secnumdepth \else
 \protect\numberline{\csname the#1\endcsname}\fi
 #7}}\fi
 \@xsect{#5}}
\def\@begintheorem#1#2{\it \trivlist \item[\hskip \labelsep{\bf #1\ #2.}]}

\def\section{\@startsection {section}{1}{\z@}{-3.5ex plus -1ex minus
 -.2ex}{2.3ex plus .2ex}{\normalsize\bf}}
\makeatother
\makeatletter
\def\subsection{\@startsection {subsection}{1}{\z@}{-3.5ex plus -1ex minus
 -.2ex}{2.3ex plus .2ex}{\normalsize\bf}}

\makeatother

\begin{document}

\begin{center}
{\large\bf A Note on Projecting the Cubic Lattice} \\
\vspace*{+.2in}

N.~J.~A.~Sloane${}^{(a)}$, Vinay~A.~Vaishampayan, \\
AT\&T Shannon Labs, \\
180 Park Avenue, Florham Park, NJ 07932-0971, USA. \\

and \\

Sueli~I.~R.~Costa, \\
University of Campinas, \\
Campinas, SP 13083-970, Brazil. \\

\vspace*{+.1in}
Email: njas@research.att.com, vinay@research.att.com, sueli@ime.unicamp.br.

\vspace*{+.1in}
${}^{(a)}$ To whom correspondence should be addressed.

\vspace*{+.1in}
April 5, 2010; revised April 24, 2010, May 25 2010, July 14 2010
\vspace*{+.1in}


{\bf Abstract}
\end{center}

It is shown that, given any $(n-1)$-dimensional lattice $\La$,
there is a vector $v \in \ZZ^n$ such that the orthogonal projection
of $\ZZ^n$ onto $v^{\perp}$ is, up to a similarity,
arbitrarily close to $\La$.
The problem arises in attempting to find the largest cylinder
anchored at two points of $\ZZ^n$ and containing no
other points of $\ZZ^n$.

\vspace{0.8\baselineskip}
Keywords: projections, shadows, dense packings

AMS 2010 Classification: Primary 11H55, 52C17


\section{Introduction}\label{Sec1}
Let $\ZZ^n$ denote the cubic lattice with basis 
$e_1 := (1,0,\ldots,0), \ldots, e_n := (0,0,\ldots,1)$.
If we project $\ZZ^n$ onto the $(n-1)$-dimensional subspace 
$$
v^{\perp} := \{ x \in \RR^n : x\cdot v = 0\}
$$ 
perpendicular to a vector $v \in \ZZ^n$,
we obtain an $(n-1)$-dimensional lattice that we denote by $\La _v$.
We will show that, given any $(n-1)$-dimensional lattice $\La$,
we can choose $v \in \ZZ^n$ so that $\La _v$ is arbitrarily
close to a lattice that is geometrically similar to $\La$.
More precisely, we will establish:

\begin{theorem}\label{Th1}
Let $\La$ be an $(n-1)$-dimensional lattice with Gram matrix $A$
$($with respect to some basis for $\RR^{n-1})$.
For any $\eps >0$,
there exist a nonzero vector $v \in \ZZ^n$,
a basis $B$ for the $(n-1)$-dimensional lattice $\La _v$
and a number $c$ such that
if $A_v$ denotes the Gram matrix of $B$, then
\beql{Eq1}
\|A - cA_v\|_{\infty} < \eps \,.
\eeq
\end{theorem}

The theorem is at first surprising, since $A$ has 
$\binom{n+1}{2}$ degrees of freedom,
whereas $v$ has only $n$ degrees of freedom (for the explanation
see the remark following the proof of Theorem \ref{Th2}).

The problem arises from a question in communication theory
(see \S\ref{struts}), which calls for projections
$\La _v$ with high packing density.
Since both the determinant and minimal norm of a lattice
are continuous functions of the entries in the Gram matrix, so is
the packing density.\footnote{The minimal norm
$\mu$ of a $d$-dimensional lattice with Gram matrix $A$
is the minimum over all $z \in \ZZ^d$, $z \ne 0$,
of the quadratic form $z A z^{\tr}$.
It is enough to consider the finite set of $z$ in some ball
around the origin. For a given $z \ne 0$, $z A z^{\tr}$
is a continuous function of the entries of $A$,
and since the minimum of a finite set of continuous functions
is continuous, $\mu$ is
a continuous function of the entries of $A$.} The theorem therefore implies
that the packing density of $\La _v$ can be made arbitrarily close
to that of $\La$.
So if we know a dense lattice in $\RR^{n-1}$,
we can find projections that converge to it in density.

\paragraph{Remark.}
We know (see for example \cite[Cor.~8]{LDLV}) that
if $\La$ is a classically integral $(n-1)$-dimensional
lattice then $\La$ can be embedded in some odd unimodular
lattice $K$ of dimension $k \le n+2$, although for
$n \ge 7$ $K$ need not be $\ZZ^k$. In any case
this does not imply that $\La$ can be recovered as a projection of $K$.

\paragraph{Notation.}
$\La^{*}$ denotes the dual lattice to $\La$,
$A^{\tr}$ is the transpose of $A$,
and $\|A\|_{\infty} = \max_{i,j} |A_{i,j}|$.
Our vectors are row vectors.
For undefined terms from lattice theory see \cite{SPLAG}.


\section{Proof of Theorem \ref{Th1}}\label{Proof}
We begin with some preliminary remarks about the
projection lattice $\La _v$ and its dual $\La _v^{*}$.
For simplicity we will only consider projections
that use vectors of the form $v=(1,v_1, v_2, \ldots, v_{n-1}) \in \ZZ^n$.
Let $\hat{v} := (v_1, v_2, \ldots, v_{n-1})$,
$M := \|v\|^2 = 1 + \sum v_i^2$.

The matrix that orthogonally projects $\RR^n$ onto $v^{\perp}$ 
is $P := I_n - \frac{1}{M} v^{\tr} v$.
As a generator matrix $G$ for $\ZZ^n$ (expressed in terms
of $e_1, \ldots, e_n$) we take $I_n$ with
its first row replaced by $v$. Let $G_v$ be
obtained by omitting the first (zero) row of $GP$.
Then $G_v$ is an $(n-1) \times n$ generator matrix for
the projection lattice $\La _v$,
and $A_v := G_v G_v^{\tr} = 
I_{n-1} - \frac{1}{M} \hat{v}^{\tr} \hat{v}$
is its Gram matrix, with $\det \La _v = \det A_v = \frac{1}{M}$. 

It is often easier to work with the dual lattice $\La _v^{*}$.
This is the intersection of $\ZZ^n$ with the subspace $v^{\perp}$,
and has generator matrix
\begin{eqnarray}
\begin{bmatrix} -v_1 & 1 & 0 & \ldots & 0 \\ 
-v_2 & 0 & 1 & \ldots & 0 \\ 
\vdots & \vdots & \vdots & \ddots & \vdots \\ 
-v_{n-1} & 0 & 0 & \ldots & 1 \\ 
\end{bmatrix} \,,
\label{Eq2}
\end{eqnarray}
Gram matrix $A_v^{*} = I_{n-1}+\hat{v}^{\tr} \hat{v}$,
and determinant $M$.

If a sequence of matrices $T_i$ converges
in the $\|~\|_\infty$ norm to a positive-definite
matrix $T$, then $T_i^{-1}$ converges to $T^{-1}$. 
So the following theorem is equivalent to Theorem \ref{Th1}.

\begin{theorem}\label{Th2}
Let $\La$ be an $(n-1)$-dimensional lattice with Gram matrix $A$
$($with respect to some basis for $\RR^{n-1})$.
For any $\eps >0$,
there exist a nonzero vector $v \in \ZZ^n$,
a basis $B$ for the $(n-1)$-dimensional lattice $\La _v^{*}$
and a number $c$ such that
if $A_v^{*}$ denotes the Gram matrix of $B$, then
\beql{Eq3}
\|A - cA_v^{*}\|_{\infty} < \eps \,.
\eeq
\end{theorem}

\noindent{\bf Proof of Theorem \ref{Th2}.}
We may write $A=L L^{\tr}$ where $L = [L_{i,j}]$
is an $(n-1) \times (n-1)$ lower triangular matrix.
For $w = 1,2,\ldots$ let us form the $(n-1) \times n$ matrix
\begin{eqnarray}
L_w & := & - \begin{bmatrix} \lfloor wL \rfloor & \pmb{0} \end{bmatrix}
+ \begin{bmatrix} \pmb{0} & I_{n-1} \end{bmatrix} \nonumber \\
& = & 
\begin{bmatrix} 
-\lfloor w L_{1,1} \rfloor & 1 & 0 & \ldots & 0 & 0 \\
-\lfloor w L_{2,1} \rfloor & -\lfloor w L_{2,2} \rfloor & 1 & \ldots & 0 & 0 \\
\vdots & \vdots & \vdots & \ddots & \vdots & \vdots \\
-\lfloor w L_{n-1,1} \rfloor & -\lfloor w L_{n-1,2} \rfloor & -\lfloor w L_{n-1,3} \rfloor & \ldots & -\lfloor w L_{n-1,n-1} \rfloor & 1 
\end{bmatrix} \,,
\label{Eq4}
\end{eqnarray}
where $\pmb{0}$ denotes a column of $n-1$ zeros.
We call $L_w$ a ``lifted'' version of $L$.

We apply elementary row operations to $L_w$ so as to put it in the form
\begin{eqnarray} 
\tilde{L}_w := \begin{bmatrix} 
-v_1 & 1 & 0 & \ldots & 0 & 0 \\
-v_2 & 0 & 1 & \ldots & 0 & 0 \\
\vdots & \vdots & \vdots & \ddots & \vdots & \vdots \\ 
-v_{n-2} & 0 & 0 & \ldots & 1 & 0 \\
-v_{n-1} & 0 & 0 & \ldots & 0 & 1
\end{bmatrix} \,,
\label{Eq5}
\end{eqnarray}
and take $v = (1,v_1, \ldots, v_{n-1})$.
Then $\La _v^{*}$ has generator matrix $\tilde{L}_w$.
But $\tilde{L}_w$ and $L_w$ generate the same lattice.
It follows that $\La _v^{*}$ has a Gram matrix
\beql{Eq6}
A_v^{*} = L_w L_w^{\tr} = w^2 A + B = w^2 \big( A+\frac{1}{w^2}B \big) \,,
\eeq
using \eqn{Eq4}, 
where the entries in $B$ are of order $O(w)$ as $w \rightarrow \infty$.
This implies \eqn{Eq3} (with $c=1/w^2)$.~~~$\bsq$

\paragraph{Remark.}
The apparent paradox mentioned in \S\ref{Sec1}
is explained by the fact that we use $\binom{n}{2}$
degrees of freedom in going from \eqn{Eq4} to \eqn{Eq5}.


\section{Examples}\label{examples}
\subsection{The lattice $2\ZZ \oplus \ZZ$}\label{2ZZ}
We start with a concrete example.
If we take $v=(1,1,0)$ then $\La _v$ has Gram matrix
$\frac{1}{2}\begin{bmatrix} 2 & 0 \\ 0 & 1 \end{bmatrix}$,
and is geometrically similar to $\sqrt{2} \ZZ \oplus \ZZ$.
Similarly $v=(1,1,1)$ produces the hexagonal (or $A_2$)
lattice, and in general $v=(1,1,\ldots,1)$ gives $A_{n-1}$. 
On the other hand, there is no $v=(1,v_1,v_2) \in \ZZ^3$ such that
$\La _v$ is geometrically similar to $2\ZZ \oplus \ZZ$
(see Proposition \ref{Prop1}).
However, we can find projections which converge to
a lattice that is geometrically similar to $2\ZZ \oplus \ZZ$. 
Since any two-dimensional lattice is geometrically
similar to its dual, we can apply Theorem \ref{Th2}
with $\La = 2\ZZ \oplus \ZZ$. Then 
$L = \begin{bmatrix} 2 & 0 \\ 0 & 1 \end{bmatrix}$,
the lifted generator matrix is
$L_w = \begin{bmatrix} -2w & 1 & 0 \\ 0 & -w & 1 \end{bmatrix}$,
$\tilde{L}_w = \begin{bmatrix} -2w & 1 & 0 \\ -2w^2 & 0 & 1 \end{bmatrix}$,
$v=(1,2w,2w^2)$,
and a Gram matrix for $\La _v^{*}$ is $I_2 + \tilde{v}^{\tr} \tilde{v}
= \begin{bmatrix} 4w^2+1 & 4w^3 \\ 4w^3 & 4w^4 +1 \end{bmatrix}$.
If we subtract $w$ times the first generator from the second,
this becomes
$$
\begin{bmatrix} 4w^2+1 & -w \\ -w & w^2+1 \end{bmatrix}
= w^2 
\begin{bmatrix} 4+1/w^2 & -1/w \\ -1/w & 1+1/w^2 \end{bmatrix} \,,
$$
which converges to $w^2 \begin{bmatrix} 4 & 0 \\ 0 & 1 \end{bmatrix}$
as $w \rightarrow \infty$.

\begin{prop}\label{Prop1}
There is no vector $v=(1,a,b) \in \ZZ^3$ such that
$\La^{*} _v$ is geometrically similar to $2\ZZ \oplus \ZZ$.
\end{prop}

\noindent{\bf Proof.}
From \eqn{Eq2}, $\La^{*} _v$ has Gram matrix 
$A := \begin{bmatrix} a^2+1 & ab \\ ab & b^2+1 \end{bmatrix}$.
If $\La^{*} _v$ is geometrically similar to $2\ZZ \oplus \ZZ$
then there is a matrix $T := \begin{bmatrix} r & s \\ t & u \end{bmatrix}
\in SL_2(\ZZ)$ and $ \lambda  \in \RR$
such that 
$$
A =  \lambda  \, T \begin{bmatrix} 4 & 0 \\ 0 & 1 \end{bmatrix} T^{\tr}
=  \lambda  \begin{bmatrix} 4r^2+s^2 & 4rt+su \\ 4rt+su & 4t^2+u^2 \end{bmatrix} \,.
$$
This implies $ \lambda  \in \QQ$, and taking the determinant and trace
of both sides we obtain $a^2+b^2+1=4 \lambda ^2$, $a^2+b^2+2=4 \lambda ^2+1= \lambda  \si$,
where $\si := 4r^2+4t^2+s^2+u^2 \in \ZZ$.
Hence the discriminant of the quadratic for $ \lambda $, $\si ^2 - 16$, is a perfect square,
so $\si = 4$ or $5$,
$ \lambda  =\frac{1}{2}$, $\frac{1}{4}$ or $1$, $a^2+b^2+1=1$, $\frac{1}{4}$ or $3$,
none of which are possible.~~~$\bsq$

\subsection{The lattice $5_1$}\label{51}
For an example where the floor operations in \eqn{Eq4}
are actually needed, consider the lattice $\La$
with Gram matrix $\begin{bmatrix} 3 & 1 \\ 1 & 2 \end{bmatrix}$
and determinant $5$ (this is the lattice $5_1$ in
the notation of \cite{LDLI}).
Again there is no $v=(1,v_1,v_2) \in \ZZ^3$ such that
$\La _v$ is geometrically similar to $\La$.
We take $L = 
\begin{bmatrix} \sqrt{3} & 0 \\ \frac{1}{\sqrt{3}} & \sqrt{\frac{5}{3}} \end{bmatrix}$,
and find that
$$
v = \big(1, \lfloor \sqrt{3} \, w \rfloor,
\lfloor \sqrt{3} \, w \rfloor \lfloor \sqrt{5/3} \, w \rfloor +
\lfloor w/\sqrt{3} \rfloor \big) \,.
$$

\subsection{The lattice $D_m$, $m \ge 3$}\label{Dn}
As generator matrix for $D_m^{*}$ we take (\cite[p.~120]{SPLAG})
\beql{EqDn}
\begin{bmatrix} 
1 & 0 & 0 & \ldots & 0 & 0\\
0 & 1 & 0 & \ldots & 0 & 0 \\
\vdots & \vdots & \vdots & \ddots & \vdots & \vdots \\
0 & 0 & 0 & \ldots & 1 & 0 \\
1/2 & 1/2 & 1/2 & \ldots & 1/2 & 1/2
\end{bmatrix} \,.
\eeq
We set $w=2t$, $t \in \ZZ$, and obtain
$$
v=\left(1,2t,(2t)^2, \ldots,(2t)^{m-1},t \frac{(2t)^m-1}{2t-1}\right) \,.
$$
In particular, when $m=3$, we have
\beql{EqFCC}
v = (1,2t,4t^2,4t^3+2t^2+t) \,,
\eeq
for which $\La _v^{*}$ converges to
the body-centered cubic lattice $D_3^{*}$
and $\La _v$ to the face-centered cubic lattice $D_3$.

\subsection{The lattice $E_8$}\label{E8}
Using the generator matrix given in \cite[p.~121]{SPLAG}, we 
find that 
$v=(1,v_1,v_2,\ldots,v_8)$ is given by
$v_1 = 2w$,
$v_2 = 2w^2-w$,
$v_i = w(v_{i-1}-v_{i-2})$ for $i=3,4,\ldots,7$,
and $v_8 = (w/2)(1 + \sum_{i=1}^7v_i)$,
where $w$ is even.

\subsection{The Leech lattice}\label{Leech}
Using \cite[Fig.~4.12]{SPLAG},
we find that
$v=(1,v_1,v_2,\ldots,v_{24})$ is given by
\begin{eqnarray}
v_1 & = & 8w, \nonumber \\
v_i & = & 4w(v_{i-1}+1), \mbox{~for~} i=2,\ldots,7,9,10,11,13,17, \nonumber \\
v_i & = & 2w(\sum_{j \in S_i}v_j +1), \mbox{~for~}i=8,12,14,15,16,18,19,20, \nonumber \\
v_i & = & 2w \sum_{j \in S_i}v_j, \mbox{~for~}i=21,22,23, \nonumber \\
v_{24} & = & w(\sum_{i=1}^{23} v_i -3),
\end{eqnarray}
where $S_8=\{1,2,\ldots,7\}$, $S_{12}=\{1,2,3,8,9,10,11\}$, $S_{14}=\{1,4,5,8,9,12,13\}$, $S_{15}=\{2k, 1 \leq k \leq 7\}$, $S_{16}=\{ 3,4,7,8,11,12,15\}$, $S_{18}=\{2,4,7,8,9,16,17\}$, $S_{19}=\{3,4,5,8,10,16,18\}$, $S_{20}=\{ 1,4,6,8,11,16,19\}$, $S_{21}=\{1,2,3,4,8,12,16,20\}$, $S_{22}=\{8,9,12,13,16,17,20,21\}$, $S_{23}=\{2k,k=4,5,\ldots,11\}$.


\section{Faster convergence}\label{faster}
The construction in Theorem \ref{Th2} produces
a vector $v$ of length $\|v\| = O(w^n)$,
while from \eqn{Eq6} we have 
$\|A - \frac{1}{w^2} A_v^*\|_{\infty} = O(\frac{1}{w}) 
= O(\frac{1}{\|v\|^{1/n}})$.
It is sometimes possible to obtain a faster rate of convergence.
Suppose $\La$ is $D_3$, and
instead of \eqn{EqDn} let us take the following
generator matrix for $D_3^{*}$:
$$
L :=
\begin{bmatrix}
-1 & -1 & ~1 \\
~1 & -1 & ~1 \\
~1 & ~1 & ~1
\end{bmatrix} \,.
$$
Let
\beql{fudge}
L_w = 
\begin{bmatrix} 
w-1 & w+1 & -w & 0 \\ 
-w-1 & w & -w+1 & 0 \\ 
-w & -w & -w & 1
\end{bmatrix}
=
\begin{bmatrix}
-wL & \pmb{0} 
\end{bmatrix} + 
\begin{bmatrix} \pmb{0} & I_3\end{bmatrix} 
+ 
\begin{bmatrix} 
-1 & 0 & 0 & 0 \\
-1 & 0 & 0 & 0 \\ 
0 & 0 & 0 & 0
\end{bmatrix} \,,
\eeq
with
\beql{fudge2}
A_v^*=L_w L_w^{\tr} = 
\begin{bmatrix}
3w^2 + 2 & w^2 + 1 & -w^2 \\ 
w^2+1 & 3w^2 + 2 & w^2 \\ 
-w^2 & w^2 & 3w^2+1 
\end{bmatrix} \,.
\eeq
The last matrix in \eqn{fudge} is chosen so that 
there are no terms of order $w$ in \eqn{fudge2}. 
Let $H_w$ denote the $3 \times 3$ matrix formed 
by the last three columns of $L_w$, and define
$v_1, v_2, v_3$ by
$$
H_w^{-1} L_w = 
\begin{bmatrix}
-v_1 & 1 & 0 & 0 \\
-v_2 & 0 & 1 & 0 \\
-v_3 & 0 & 0 & 1 
\end{bmatrix} \,.
$$
Then
\beql{vfcc}
v = (1,v_1,v_2,v_3) = (1,2w^2-w+1,2w^2+w+1,4w^3+3w)
\eeq
has $\|v\|=O(w^3)$,
and now $\|A - \frac{1}{w^2} A_v^*\|_{\infty} = O(\frac{1}{\|v\|^{2/3}})$,
which is a faster convergence than we found in \S\ref{Dn}.
We do not know if similar speedups are always possible. 
Incidentally, we first found \eqn{vfcc}---before
Theorem \ref{Th1} was proved---by a combination
of computer search and guesswork.


\section{The fat strut problem}\label{struts}
The problem studied in this paper arose when
constructing codes for a certain analog
communication channel \cite{VaishampayanCosta:2003}.
The codes require that one finds a curved tube in the sphere $S^{2n-1}$
which does not intersect itself, 
has a specified length and as large a volume as possible.
The method used in \cite{VaishampayanCosta:2003}
is based on finding a vector $v \in \ZZ^n$ with a specified value of $\|v\|$,
such that there is a cylinder of large volume with axis 
$\overrightarrow{0v}$ which contains no points of 
$\ZZ^n$ other that $0$ and $v$.
The cross-section of the cylinder is an $(n-1)$-dimensional ball,
and $0$ and $v$ are the centers of the two end-faces.
The radius of the cylinder is taken to be as large as possible
subject to the condition that the interior contains
no point of $\ZZ^n$.
The problem is to choose $v$, for a given length $\|v\|$,
so that the volume of the resulting cylinder is maximized.
We call a cylinder which achieves the
maximal volume a {\em fat strut}.

A fat strut has the property that the projection of
the cylinder onto $v^{\perp}$ does not contain the image
of any nonzero point of $\ZZ^n$.
The radius of the cylinder is therefore equal to
the radius of the largest $(n-1)$-dimensional sphere 
in the projection lattice $\La _v$ which contains 
no nonzero point of $\La _v$.
In other words, finding a fat strut for a given length $\|v\|$ 
is equivalent to maximizing the density of the 
projection lattice $\La _v$.

It is worth contrasting the fat strut problem with the result of
\cite{Heppes:1960} and \cite{Horvath:1970} that for any lattice sphere packing
in dimension three or higher there is always an {\em infinite} cylinder
of nonzero radius (obviously not passing through the origin)
which does not touch any of the spheres.

\section{Acknowledgment}
We thank the referee for some very helpful comments.

\end{document}